\documentclass[12pt]{article}
\usepackage[centertags]{amsmath}
\usepackage{amsfonts}
\usepackage{amssymb}
\usepackage{amsthm}
\usepackage{amsmath,mathrsfs}
\usepackage{amsmath,amscd}
\usepackage[matrix,arrow,curve,cmtip]{xy}
\title{\textbf{$\infty$-Topoi and Natural Phenomena: Generation}}
\author{Renaud Gauthier \footnote{rg.mathematics@gmail.com} \\ \\}
\theoremstyle{definition}

\newcommand{\beq}{\begin{equation}}
\newcommand{\eeq}{\end{equation}}
\newcommand{\rarr}{\rightarrow}

\newcommand{\cB}{\mathcal{B}}
\newcommand{\cC}{\mathcal{C}}
\newcommand{\cD}{\mathcal{D}}
\newcommand{\cE}{\mathcal{E}}

\newcommand{\cH}{\mathcal{H}}

\newcommand{\cO}{\mathcal{O}}
\newcommand{\cP}{\mathcal{P}}
\newcommand{\cR}{\mathcal{R}}
\newcommand{\cS}{\mathcal{S}}
\newcommand{\cT}{\mathcal{T}}

\newcommand{\cU}{\mathcal{U}}

\newcommand{\cX}{\mathcal{X}}
\newcommand{\cY}{\mathcal{Y}}
\newcommand{\gC}{\mathfrak{C}}

\newcommand{\Cat}{\text{Cat}}
\newcommand{\Catinf}{\Cat_{\infty}}
\newcommand{\CatD}{\Cat_{\Delta}}
\newcommand{\diag}{\text{diag}}
\newcommand{\Dop}{\Delta^{\op}}

\newcommand{\ffqc}{\text{ffqc}}
\newcommand{\Fun}{\text{Fun}}
\newcommand{\Hom}{\text{Hom}}

\newcommand{\Ho}{\text{Ho}}
\newcommand{\inj}{\text{inj}}
\newcommand{\Loc}{\text{Loc}}
\newcommand{\loc}{\text{loc}}
\newcommand{\Kan}{\mathcal{K}an}
\newcommand{\Map}{\text{Map}}
\newcommand{\mono}{\text{mono}}

\newcommand{\op}{\text{op}}

\newcommand{\SeT}{\text{SeT}}
\newcommand{\Set}{\text{Set}}

\newcommand{\Sh}{\text{Sh}}
\newcommand{\Top}{\text{Top}}
\newcommand{\uHom}{\underline{\Hom}}

\newcommand{\dkAff}{\text{d}k\text{-Aff}}
\newcommand{\dStk}{\text{dSt}(k)}
\newcommand{\dStkinf}{\dStk_{\infty}}

\newcommand{\LBous}{L_{\text{Bous}}}
\newcommand{\LDK}{L_{\text{DK}}}
\newcommand{\OXmono}{\cO_{\cX}^{\text{(mono)}}}

\newcommand{\RHom}{\mathbb{R} \uHom}

\newcommand{\RHomg}{\RHom_{\SeT}}

\newcommand{\RHomLg}{\RHom^*_{\SeT}}
\newcommand{\SetD}{\Set_{\Delta}}
\newcommand{\SeCat}{\text{SeCat}}
\newcommand{\SePC}{\text{SePC}}

\newcommand{\skCAlg}{\text{s}k\text{-CAlg}}

\newcommand{\sPr}{\text{sPr}}
\newcommand{\sPrtau}{\text{sPr}_{\tau}}

\begin{document}
\maketitle
\begin{abstract}
We show that the Segal topos being used in \cite{RG} to model natural phenomena has a subobject classifier, something we regard as being a source from which dynamics is generated. This is done by considering the $\infty$-category associated to such a Segal topos in the language of \cite{Lu1}, which turns out to be an $\infty$-topos. At this point we have the formalism of Higher topoi at our disposal to deal with Higher Category Theory concepts in a transparent manner.
\end{abstract}

\newpage

\section{Introduction}
In \cite{RG} we argued that one could model natural laws by simplicial commutative algebras and that representations of such laws could be implemented by considering higher stacks, or derived stacks (\cite{TV4}, \cite{T}) to be more precise. For $k$ a commutative ring, $\skCAlg$ the simplicial category of simplicial commutative $k$-algebras, which after Segal localization (\cite{T}, \cite{TV1}) gives $\dkAff = (L(\skCAlg))^{\op}$ the Segal category of derived affine schemes, one can first consider the Segal category of derived pre-stacks on $\dkAff$, $\widehat{\dkAff} = \RHom(\dkAff^{\op}, \Top)$. On $\dkAff$ we put the ffqc topology, and we define the Segal category of derived stacks $\dStk$ on $\dkAff$, left Bousfield localization of $\widehat{\dkAff}$ with regards to ffqc-hypercovers. This is a Segal topos (\cite{T}). The purpose of this paper is to discuss the existence of subobject classifiers in such a Segal topos. It is in the setting of $\infty$-topoi that such a concept is easier to discuss, so we will consider the $\infty$-topos associated to $\dStk$. That we have such an association is a fairly well-known fact; as pointed out in \cite{Lu2}, \cite{Lu1}, \cite{TV3}, \cite{TV4}, the notions of Segal topos and $\infty$-topos are essentially equivalent. To be precise, we have an adjunction from simplicial categories to simplicial sets, with right adjoint the simplicial nerve functor $N$ (\cite{Lu1}), originally introduced by Cordier as the coherent nerve functor (\cite{C}). We have $N(\dStk) \in \Catinf$. We show $\dStkinf = N(\dStk)$ is an $\infty$-topos. Focusing then on subobject classifiers in $\infty$-topos, we first argue that an object of an $\infty$-topos, viewed as a model for a realization of natural laws, is fully defined by its subobjects. From this perspective, one regards subobject classifiers as generators of natural phenomena in a given $\infty$-topos, and since $\infty$-topos have such an object, we conclude that $\dStkinf$, the $\infty$-category associated to $\dStk$, has a subobject classifier, hence a source from which phenomena are drawn. \\

In \cite{RG} the emphasis was on Higher Galois, and for this purpose we considered $\RHomLg(\cX, \cU)$ for two Segal topos $\cX$ and $\cU$, and we showed those are Segal groupoids. From the perspective of the present paper it is preferable to have Segal topos however, and $\RHom(\cX, \cU)$ would be a more apt object. In the context of $\infty$-topos we show $\Fun(\cX, \cU)$ is an $\infty$-topos for two $\infty$-topoi $\cX$ and $\cU$. In particular so is $\Fun(\cX,\cX)$, so it has a subobject classifier. The conclusions we drew in \cite{RG} therefore hold on the $\infty$-topoi side of the problem as well. From there one investigates whether there is a universal such $\infty$-topos, universal in the sense of having a universal subobject classifier, and we show that indeed there is one, namely the $\infty$-category of spaces $\cS$.\\

In the present paper, whatever technical definition is invoked, but is not fully used, will be referenced, and only those concepts that need to be explained will be covered, for the sake of ease of reading. In a first time we go over the necessary background about Segal categories and $\infty$-categories. This we do in Section 2. As an intermission we discuss localizations (Dwyer-Kan, Segal, Hammock, Bousfield) in Section 3, at which point we can come back to Segal categories and $\infty$-categories to define those localizations, which leads to the notions of Segal topoi and $\infty$-topoi. We then define the Segal categories of pre-stacks and stacks with this formalism. Moving on to the $\infty$-category side of things in Section 4, we then use a result of \cite{Lu1} which essentially states that the simplicial nerve of the subcategory of fibrant-cofibrant objects in the model category of prestacks in the local model structure, meaning the Segal category of stacks, is an $\infty$-topos, which is what we needed to show that $\dStk$ maps to an $\infty$-topos. From there we get results mirroring those of \cite{RG}, and we then discuss subobject classifiers and how relevant $\infty$-topoi relate to one another in this regard.\\

For notations, $\SetD$ is the category of simplicial sets, $\Delta$ is the usual category of combinatorial simplices whose objects are linearly ordered sets $[n] = \{0,1, \cdots, n \}$, with morphisms order preserving maps. For $A: \Delta \rarr \Set$, we write $A_n$ for $A([n])$. Our main reference for model categories will be \cite{Ho}.\\

\section{Higher categories}
\subsection{Segal Categories}
Our main references for Segal categories will be \cite{HS}, \cite{P}, \cite{TV1}. One can regard a Segal category as a weak form of simplicial category, where compositions are defined up to equivalence only. To be more precise, a Segal category $C$ is first a bi-simplicial set, a functor $C : \Dop \rarr \SetD$, such that $C_0$ is a discrete set, its set of objects. This makes a bisimplicial set into a Segal pre-category. The simplicial set of morphism $C_1(x,y)$ between objects $x,y \in C_0$ is the pre-image of $(x,y) \in C_0 \times C_0$ under the morphism $C_1 \rarr C_0 \times C_0$. What characterizes a Segal category among Segal pre-categories is the fact that the morphisms of simplicial sets $C_n \rarr C_1 \times_{C_0} \cdots \times _{C_0} C_1$ are equivalences for all $n \in \Delta$. It is in this regard that Segal categories are generalizations of simplicial categories. Indeed:
\beq
C_n \cong \coprod_{(a_0, \cdots, a_n) \in C_0^{n+1}} C_{(a_0, \cdots, a_n)} \nonumber
\eeq
so for all points $(a_0, \cdots, a_n) \in C_0^{n+1}$, a Segal category would just provide us with equivalences of simplicial sets:
\beq
C_{(a_0, \cdots, a_n)} \rarr C_{(a_0, a_1)} \times \cdots \times C_{(a_{n-1}, a_n)} \nonumber
\eeq
whereas if say $S$ is a simplicial category (\cite{GJ}), defining:
\beq
S_n = \coprod_{(x_0, \cdots, x_n) \in S_0^{n+1}} \uHom(x_0,x_1) \times \cdots \times \uHom(x_{n-1}, x_n) \nonumber
\eeq
we have natural isomorphisms of simplicial sets:
\beq
S_{(x_0, \cdots, x_n)} \rarr S_{(x_0, x_1)} \times \cdots \times S_{(x_{n-1}, x_n)} \nonumber
\eeq
To come back to Segal pre-categories, they form a category $\SePC$ on which we can put two model categories, and for this we refer the reader to \cite{TV1}. One of those is a left Bousfield localization of the other, and the natural inclusion on the corresponding homotopy categories has a left adjoint $\SeCat$ that turns Segal pre-categories into Segal categories. The main point here is that $\SeCat$ is a left derived functor, a fact we will use later. We need a notion of homotopy category $\Ho(C)$ for a Segal category $C$: it suffices to take a category with $C_0$ as set of objects and $\pi_0(C_{(a,b)})$ as set of morphisms between two objects $a,b \in C_0$. A morphism $f: C \rarr D$ of Segal pre-categories will be said to be fully faithful if for all $a,b \in C_0$ the morphism $\SeCat(C)_{(a,b)} \rarr \SeCat(D)_{(f(a),f(b))}$ is an isomorphism in $\Ho(\SetD)$, and $f$ is said to be essentially surjective if $\Ho(\SeCat(C)) \rarr \Ho(\SeCat(D))$ is essentially surjective. As usual, we will define a morphism of Segal pre-categories to be an equivalence if it is both fully faithful and essentially surjective. If we take this as a notion of weak equivalence on $\SePC$ and we take cofibrations to be monomorphisms, this defines a model structure on $\SePC$ as shown by Pellissier (\cite{P}, \cite{TV1}). This model structure is furthermore internal, with an internal hom $\uHom$, and $\Ho(\SePC)$ is cartesian closed, with internal hom $\RHom$. When we say $f: C \rarr D$ is a morphism of Segal categories, we will mean it's an element of $\RHom(C,D)$. Very useful also is the derived adjunction formula:
\beq
\RHom(A, \RHom(B,C)) \cong \RHom(A \times B, C) \nonumber
\eeq

\subsection{$\infty$-categories}
A standard reference for our purposes will be \cite{Lu1}. By $\infty$-category, we mean $(\infty,1)$-category, where a $(\infty,n)$-category is an $\infty$-category for which morphisms are invertible for $k > n$. We say a simplicial set $K$ is an $\infty$-category if one can find a dotted arrow as shown below that can make such a diagram commutative for $0 < i < n$, where $\Delta^n = \Hom_{\SetD}(-,[n])$ and $\Lambda_i^n$ is the $i$-th horn of $\Delta^n$:
\beq
\xymatrix{
\Lambda_i^n \ar@{^{(}->}[d] \ar@{->}[r] &K \\
\Delta^n \ar@{.>}[ur] \nonumber
}
\eeq
Those simplicial sets which satisfy such an extension property for outer horns as well are referred to as Kan complexes. If $\cC$ and $\cD$ are $\infty$-categories, a functor from $\cC$ to $\cD$ is simply a map $\cC \rarr \cD$ of simplicial sets. We denote $\Map_{\SetD}(\cC, \cD)$ by $\Fun(\cC, \cD)$, which is the $\infty$-category of functors from $\cC$ to $\cD$. We denote by $\Catinf$ the $\infty$-category of small $\infty$-categories.\\

\subsection{Simplicial categories}
$\SetD$-enriched categories are referred to as simplicial categories, the category of which will be denoted by $\CatD$ (see \cite{GJ} as well). Crucial for us will be the simplicial nerve functor $N$:
\beq
N: \CatD \rarr \SetD \nonumber
\eeq
as originally introduced by Cordier under the name of coherent nerve functor (\cite{C}), and whose definition is given again in \cite{Lu1} in a format suitable for our purposes. It is defined via the simplicial category functor $\gC$ as presented in \cite{Lu1}. We will not need its full definition, and the reader is referred to that reference for details. One can see $\gC: \SetD \rarr \CatD$ as a thickening of simplicial sets, through a combinatorial definition. If $\cC$ is a simplicial category, the simplicial nerve $N(\cC)$ is defined by the adjunction formula:
\beq
\Hom_{\SetD}(\Delta^n, N(\cC)) = \Hom_{\CatD}(\gC[\Delta^n], \cC) \nonumber
\eeq
which one may want to contrast with the usual nerve functor with satisfies, for $\cC$ an ordinary category:
\beq
\Hom_{\SetD}(\Delta^n, N(\cC)) = \Hom_{\Cat}([n], \cC) \nonumber
\eeq
Note that we have an equivalence of simplicial categories $\gC[\Delta^n] \rarr [n]$. $\Kan$, the full subcategory of $\SetD$ spanned by Kan complexes, is a simplicial category. $\cS = N(\Kan)$ is referred to as the $\infty$-category of spaces.

\section{Localizations}
We will use a few different localizations, so for the sake of fixing ideas we will go over those we need and provide references for more ample details.
\subsection{Bousfield localization}
Our main reference for Bousfield localizations will be \cite{Hi}. Recall that if $M$ is a model category, $S$ a class of maps in $M$, an object $X$ of $M$ is said to be $S$-local, if it is fibrant and if for any $f:A \rarr B$ in $S$, the induced map of homotopy function complexes $f^*: \Map(B,X) \rarr \Map(A,X)$ is a weak equivalence in $\SetD$. Now a map $g: X \rarr Y$ in $M$ is a $S$-local equivalence if for all $S$-local object $Z$, the induced map $g^*: \Map(Y,Z) \rarr \Map(X,Z)$ is a weak equivalence in $\SetD$. We then define the left Bousfield localization $L_S M$ of $M$ to be a model structure on the underlying category of $M$ with $S$-local equivalences as weak equivalences, the same cofibrations as those of $M$, and fibrations will be maps with the right lifting property with respect to cofibrations that are also $S$-local equivalences. We will denote Bousfield localizations by $\LBous$, the set $S$ being implied.

\subsection{Simplicial localizations}
In \cite{DK1} Dwyer and Kan introduced a slightly more general localization than the usual Gabriel-Zisman localization (\cite{GZ}); for $\cC$ an ordinary category, $\cB$ a subcategory, the Gabriel-Zisman localization of $\cC$ with respect to $\cB$ is denoted $\cC[\cB^{-1}]$ (or sometimes $\cB^{-1} \cC$), obtained by formally inverting morphisms in $\cB$. Dwyer and Kan generalized this by defining the simplicial localization of $\cC$ with respect to $\cB$ to be given by $\LDK(\cC, \cB) = F_*\cC [F_* \cB^{-1}]$ where $F_*\cC$ denotes the simplicial free resolution of the category $\cC$ (see \cite{DK1}). $\LDK(\cC, \cB)$ is a simplicial category, and it is a generalization in the sense that $\Ho(\LDK[\cC, \cB]) = \cC[\cB^{-1}]$. One can also generalize this definition to the case where $\cC$ is a simplicial category by taking the diagonal of the bi-simplicial set constructed using the above procedure levelwise: $\LDK(\cC, \cB) = \text{diag} F_*\cC [F_* \cB^{-1}]$.\\

In \cite{TV1}, \cite{T}, a still more general localization is defined, using the fact that we have a fully faithful functor $\CatD \rarr \SeCat$. Toen and Vezzosi define a Segal analog of $\LDK$ as follows: if $C$ is a Segal category, $S$ a set of morphisms in $\Ho(C)$, there exists a Segal category $L(C,S)$, with a localization $l: C \rarr L(C,S)$, characterized by the following universal property: for any Segal category $D$, the induced morphism:
\beq
l^*: \RHom(L(C,S),D) \rarr \RHom(C,D) \nonumber
\eeq
is fully faithful and has for essential image those morphisms from $C$ to $D$ that send morphisms in $S$ into equivalences in $D$. Observe that if $C$ is an ordinary category viewed as a Segal category, then $L(C,S) = \LDK(C,S)$ (\cite{HS}). If $M$ is a model category, $W$ its set of equivalences, we will just denote $L(M,W)$ in the sense above by $LM$. We will denote both Segal localizations and Dwyer-Kan localizations by the same letter $L$, and it should be clear from the context which one is being used.\\

In \cite{DK2}, yet another simplicial localization is introduced, the hammock localization $L^H$, which is better behaved than the original Dwyer-Kan localization $\LDK = L$. We won't need its explicit definition so the reader is referred to \cite{DK2} for details, but what is important for us is that if $C$ denotes an ordinary category, we have a weak equivalence:
\beq
L^H C \rarr LC \nonumber
\eeq
For $C$ a simplicial category, we can also define its hammock localization as $L^H C = \diag L^H C$, and we also have a weak equivalence in this case:
\beq
\diag L^H C \rarr LC \nonumber
\eeq
From \cite{DK3}, if $M$ is a simplicial model category, $M^{\circ}$ its subcategory of fibrant-cofibrant objects, we have a weak equivalence:
\beq
M^{\circ} \rarr \diag L^H M \nonumber
\eeq
Collecting things, if $M$ is a simplicial model category, $LM \simeq \diag L^H M \simeq M^{\circ}$ so $LM \simeq M^{\circ}$ as observed in \cite{HS}. This is a crucial point that will be important later.

\subsection{Higher topoi}
\subsubsection{Segal topoi}
As defined in \cite{TV1}, a map of Segal categories is said to be left exact if it preserves finite limits. A localization of Segal categories is a morphism of Segal categories with a fully faithful right adjoint. A Segal category $\cX$ is a Segal topos if it is a left exact localization of $\hat{C} = \RHom(C^{\op} , \Top)$, with $\Top = L\SetD$, the Segal category of simplicial sets, and $C$ a small Segal category. This means there exists some $i: \cX \rarr \hat{C}$ fully faithful with a left exact left adjoint. As in classical topos theory (\cite{MM}), Toen and Vezzosi define morphisms between Segal topoi to be geometric morphisms; if $\cX$ and $\cY$ are Segal topoi, one denotes by $\RHomg(\cX, \cY) = \RHomLg(\cY, \cX)$ the Segal category of geometric morphisms from $\cX$ to $\cY$, defined as the sub-Segal category of $\RHom(\cY, \cX)$ spanned by those morphisms which are left exact with a right adjoint. We denote the 2-Segal category of Segal topoi by $\SeT$.

\subsubsection{$\infty$-topoi}
As defined in \cite{Lu1}, a localization is a functor between $\infty$-categories that has a fully faithful right adjoint. The notion of adjoints for $\infty$-categories makes use of a correspondence, and is an abstraction of the notion of adjunction from ordinary category theory (\cite{McL}). Nevertheless, this definition, which we won't use, can be recast in a classical format by virtue of the fact that if $f: \cC \rarr \cD$ and $g: \cD \rarr \cC$ are functors between $\infty$-categories $\cC$ and $\cD$, $f \dashv g$ if and only if there exists a morphism $\eta: id_{\cC} \rarr g \circ f$ in $\Fun(\cC, \cC)$ such that for all $X \in \cC$, $Y \in \cD$, we have an induced equivalence in the homotopy category of spaces $\cH$:
\beq
\Map_{\cD}(f(X),Y) \simeq \Map_{\cC}(X, g(Y)) \nonumber
\eeq
If $f: \cC \rarr \cD$ is a localization, $\cD$ is said to be a localization of $\cC$. If $g$ denotes the right adjoint that makes $f$ a localization, $L = g \circ f : \cC \rarr \cC$ will be referred to as a localization functor. For the definition of right (left) exactness, the reader is referred to this reference, but the definition is not very illuminating in our context. What is important is the following fact: a functor is left exact if and only if it preserves finite limits. $\cX \in \Catinf$ is an $\infty$-topos if there exists a small $\infty$-category $\cC$ and an accessible left exact localization functor $\cP(\cC) \rarr \cX$ where $\cP(\cC) = \Fun(\cC^{\op}, \cS)$. We will not need to check that a localization is accessible for the simple reason that most of the categories we will deal with are accessible. It is a fact (\cite{Lu1}, \cite{Lu2}) that a morphism between accessible $\infty$-categories with an adjoint is accessible, so we will dispense with this accessibility condition. Nevertheless we can briefly remind the reader of what it means to be accessible (see \cite{Lu1} for more details). We first have to go back to simplicial sets. We define the join $K \star K'$ of two simplicial sets $K$ and $K'$ as follows: for every nonempty, finite, linearly ordered set $I$:
\beq
(K \star K')(I) = \coprod_{I = J \cup J'} K(J) \times K'(J') \nonumber
\eeq
With this in hand, we denote $K \star \Delta^0$ by $K^{\triangleright}$ and refer to it as the right cone of $K$. Now if $\kappa$ is a regular cardinal, $\cC \in \Catinf$, we say $\cC$ is $\kappa$-filtered if for any small $K \in \SetD$, for any map $f: K \rarr \cC$, there is an extension $\overline{f}: K^{\triangleright} \rarr \cC$ of $f$. With this we have an obvious notion of $\kappa$-filtered colimits. For $\cC \in \Catinf$ with small $\kappa$-filtered colimits then, we say a functor $f: \cC \rarr \cD$ in $\Catinf$ is $\kappa$-continuous if it preserves $\kappa$-filtered colimits. For $\cC$ a small $\infty$-category, we let $\text{Ind}_{\kappa}(\cC)$ be the full subcategory of $\cP(\cC)$ spanned by functors $f: \cC^{\op} \rarr \cS$ that classify right fibrations $\tilde{\cC} \rarr \cC$ where $\tilde{\cC}$ is a $\kappa$-filtered $\infty$-category. Finally $\cC \in \Catinf$ is said to be $\kappa$-accessible, or simply accessible, if there exists a small $\infty$-category $\cC_0$ and an equivalence $\text{Ind}_{\kappa}(\cC_0) \rarr \cC$. If $\cC$ is such a category, a functor $f: \cC \rarr \cD$ is accessible if $\kappa$-continuous for some $\kappa$, a regular cardinal. Observe that if $K \in \SetD$, $\cC$ is an accessible $\infty$-category, then $\Fun(K, \cC)$ is accessible as well. That gives us a lot of accessible categories.

\subsection{local model structures}
One model structure we will put on our simplicial presheaves is the local model structure, something originally introduced by Jardine (\cite{J}), and used in \cite{TV3} and \cite{TV4} to construct Segal categories of derived stacks. First it is important to go briefly over those constructions for the sake of putting them in perspective, referring the reader to the main reference \cite{TV3} for more details. In that reference, simplicial categories are referred to as $S$-categories. Let $T$ be such a category. A $S$-topology $\tau$ on $T$ is nothing but a Grothendieck topology $\tau$ on its homotopy category $\Ho(T)$. This makes $(T, \tau)$ into what is called an $S$-site. On the category of simplicial presheaves on $T$, $\sPr(T)$, there exists a model structure, the local model structure, where equivalences are local equivalences, also known as $\pi_*$-equivalences. We denote by $\sPrtau(T)$ this new model structure, which turns out to be the left Bousfield localization of $\sPr(T)$ along local equivalences. $\sPrtau(T)$ is called the model category of stacks on $(T,\tau)$. Actually, since $id: \sPr(T) \rarr \sPrtau(T)$ preserves homotopy fiber products, we call such a localization a left exact Bousfield localization of $\sPr(T)$. We define a model topos to be a model category that is Quillen equivalent to a left exact Bousfield localization of $\sPr(T)$, $T$ a simplicial category. In particular that would make $\sPrtau(T)$ a model topos. We have the following result (\cite{TV1}): if $M$ is a model topos, then $LM$ is a Segal topos. Thus, $L(\sPrtau(T))$ is a Segal topos.\\

Independently, we can develop the same formalism of stacks starting not from a simplicial category, but from a Segal category. One defines a Segal topology on a Segal category $\cX$ to be a Grothendieck topology on its homotopy category $\Ho(\cX)$. One defines a notion of local equivalence on $\hat{\cX} = \RHom(\cX^{\op}, \Top)$ relative to the topology $\tau$ on the Segal site $(\cX, \tau)$. A stack $F$ on $\cX$ is an element of $\hat{\cX}$ such that for any local equivalence $\alpha: G \rarr H$ in $\hat{\cX}$, the induced morphism $\alpha^*: \Hom_{\Ho(\hat{\cX})}(H,F) \rarr \Hom_{\Ho(\hat{\cX})}(G,F) $ is bijective. We denote by $\cX^{\sim,\tau}$ the Segal category of stacks on $\cX$. Now it turns out the inclusion $\cX^{\sim, \tau} \rarr \hat{\cX}$ is a left exact left adjoint, making $\cX^{\sim,\tau}$ into a Segal topos.\\

Now there is a very important result in the theory of Segal categories called strictification (\cite{HS}, \cite{TV1}, \cite{TV3}, \cite{TV4}) that states that for $T$ a small simplicial category, $M$ a cofibrantly generated, simplicial model category, if we endow $M^T$ with its projective model structure, then we have an isomorphism in $\Ho(\SePC)$
\beq
L(M^T) \cong \RHom(T, LM) \nonumber
\eeq
This is sometimes written in the following format(\cite{TV3}): for $C$ a category with a subset of morphisms $S$:
\beq
L(M^{C,S}) \cong \RHom(L(C,S), LM) \nonumber
\eeq
where $M^{C,S}$ is the model category of restricted diagrams, a left Bousfield localization of $M^C$ with respect to morphisms based in $S$. In this isomorphism, we used the Quillen equivalence $M^{C,S} \simeq M^{L(C,S)}$.\\

Using this one can show:
\beq
L(\sPr(T)) \cong \RHom(T^{\op}, \Top) \nonumber
\eeq
as well as:
\beq
L(\sPrtau(T)) \cong T^{\sim,\tau} \nonumber
\eeq
Modulo isomorphisms then, we have the following commutative diagram that sums up this picture:
\beq
\xymatrix{
\sPr(T) \ar[r]^{\LBous} \ar[d]_L &\sPrtau(T) \ar[d]^L\\
\hat{T} = \RHom(T^{\op}, \Top) \ar[r]_-{\LBous} & T^{\sim,\tau}
} \nonumber
\eeq
For instance:
\begin{align}
L(\sPr((\skCAlg)^{\op})) &= L(\SetD^{\skCAlg}) \nonumber \\
&\simeq \RHom(L(\skCAlg), L(\SetD)) \nonumber \\
&= \RHom(\dkAff^{\op}, \Top) \nonumber \\
&= \widehat{\dkAff} \nonumber
\end{align}

\subsection{Segal category of derived stacks}
As covered in \cite{TV2}, \cite{TV3}, \cite{T}, for $k$ a commutative ring, $\skCAlg$ the category of simplicial commutative $k$-algebras, $\dkAff = L(\skCAlg)^{\op}$ the Segal category of derived affine schemes, on which we put the ffqc topology, $\dkAff^{\sim, ffqc} = \dStk $ is a left exact localization of $\widehat{\dkAff} = \RHom( \dkAff^{\op}, \Top)$. This will be the primary object of study for us, the reason being, as briefly discussed in \cite{RG}, that one can model natural laws by simplicial algebras, the realizations of which can, for the sake of coherence, be represented as stacks. The collection of all such representations is $\cX = \dStk$, and would correspond to natural phenomena, globally speaking.

\section{From Segal topoi to $\infty$-topoi}
It is clear that the transition from Segal topoi to $\infty$-topoi will be made using the simplicial nerve functor $N$. In particular we will be using the following result from \cite{Lu1} and \cite{Lu2}: if $\cC$ is a small category with a Grothendieck topology, $A = \sPr(C)$ the category of simplicial presheaves over $\cC$ with its local model structure (\cite{J}), $A^{\circ}$ the subcategory of fibrant-cofibrant objects of $A$ for that local model structure, then $N(A^{\circ}) \simeq \Sh(N(\cC))^{\wedge}$ where $\Sh(\cD) \subseteq \cP(\cD)$ for $\cD$ a small $\infty$-category with a Grothendieck topology (\cite{Lu1}) corresponds to a notion of higher category of sheaves (\cite{Lu1}). It turns out $\Sh(\cD)$ is a topological localization of $\cP(\cD)$ (\cite{Lu1}) and those are accessible and left exact, hence $\Sh(\cD)$ is an $\infty$-topos itself. If we further take the Bousfield localization of such an $\infty$-topos with respect to $\infty$-connective morphisms (\cite{Lu1}), one gets the hypercompletion $\Sh(\cD)^{\wedge}$, which again is an $\infty$-topos. This result above starting from an ordinary category holds also in the case where $\cC$ is a simplicial model category, and that will be important for us since we are considering taking $\cC = \skCAlg$.\\

The situation can be summarized as follows. If one takes $\cC = (\skCAlg)^{\op}$, then $A = \sPr(\cC^{\op})$ is a simplicial model category, on which we can put a local model structure. Actually on $A$ itself we first put the injective model structure, which is Quillen equivalent to the projective model structure. Then we have:
\beq
\xymatrix{
A_{\inj} \ar[d]_{N} \ar[rr]^{\LBous} && A_{\inj, \loc} \supset A^{\circ} \ar[d]^N \\
N(A_{inj}) \ar[d]|-{\simeq} \ar[rr]^{\Loc} && N(A^{\circ}) \simeq \Sh(N(\cC))^{\wedge} \\
\cP(N(\cC)) = \Fun(N(\cC), \cS) \ar[rd]_{\text{top. loc }} \\
&\Sh(N(\cC)) \ar[uur]_{\text{acc. lex loc}}} \nonumber
\eeq

We apply this to the case $\cC = (\skCAlg)^{\op}$, $A = \sPr(\cC) = \SetD^{\skCAlg}$ so that $A_{\text{inj,loc}} = \sPr_{\ffqc}(\cC)$. Now:
\beq
A^{\circ} \simeq LA = L(\sPr_{ffqc}((\skCAlg))^{\op}) \simeq \dStk \nonumber
\eeq
so that:
\beq
N(A^{\circ}) \simeq N(\dStk) = \dStkinf \nonumber
\eeq
is an $\infty$-topos since it is equivalent to an $\infty$-topos $\Sh(N(\dkAff))$. From then on we will be working with $\dStkinf$, the $\infty$-topos associated with the Segal topos $\dStk$.

\section{Subobject classifiers}
First a bit of definitions from \cite{Lu1}. If $\cX$ is an $\infty$-category with pullbacks, $S$ a class of morphisms in $\cX$, $\cO_{\cX} = \Fun(\Delta^1, \cX)$, then we denote by $\cO_{\cX}^{(S)}$ the subcategory of $\cO_{\cX}$ whose objects are elements of $S$ and morphisms are pullback diagrams. For $\cX \in \Catinf$ with pullbacks, $S$ a collection of morphisms in $\cX$ stable under pullback, we say a morphism $f: X \rarr Y$ classifies $S$ if it is a final object of $\cO_{\cX}^{(S)}$, which can equivalently be stated by saying that the object $Y$ classifies $S$. If for $S$ we take the class of monomorphisms in $\cX$, i.e. those morphisms $W \rarr Z$ that are $(-1)$-truncated objects of $\cX_{/Z}$, then a subobject classifier for an $\infty$-category $\cX$ with pullbacks is an object $\Omega$ that classifies the collection of all monomorphisms in $\cX$, or equivalently, a final object of $\OXmono$. Now every $\infty$-topos has a subobject classifier, hence so does $\dStkinf$.\\

The importance of such a concept can be traced back to the original motivation for choosing $\dStk$ in the first place; stacks model realizations of natural laws, i.e. natural phenomena. Monomorphisms into a given stack can be seen as constituents of such a stack, hence $F \in \dStk$ can be seen as being characterized by monomorphisms in $\dStk/F$. If one has a subobject classifier $\Omega$ on the $\infty$-topos side of thing, by its very definition it can be regarded as an object out of which other objects are being constructed, hence this would define a generation.\\

In \cite{RG} we worked, for $T$ and $\cX$ two Segal topos, with $\RHomLg(T,\cX)$ and showed this is a Segal groupoid only, and this was done for the sake of Higher Galois Theory. One could show that $\RHom(\cX,T)$ is a Segal topos, but since we work with $\infty$-topos here, we can just consider $\Fun(\cX, \cU)$ for $\cX, \cU $ $\infty$-topoi. In particular $\cU$ is a left exact localization of $\cP(\cC)$ for some small $\infty$-category $\cC$. From \cite{Lu2} we know $\cP(\cC)$ is an $\infty$-topos. Denote by $\pi: \cP(\cC) \rarr \cU$ the left exact left adjoint of this localization. It induces a localization $\pi_*: \Fun(\cX, \cP(\cC)) \rightleftarrows \Fun(\cX, \cU):i_*$, again left exact. Now:
\begin{align}
\Fun(\cX, \cP(\cC)) &= \Fun(\cX, \Fun(\cC^{\op}, \cS)) \nonumber \\
&\simeq \Fun(\cX \times \cC^{op}, \cS) \nonumber
\end{align}
and the latter is a $\infty$-topos, hence so is $\Fun(\cX, \cU)$ by localization, and in particular so is $\Fun(\cX, \cX)$, where for us $\cX = \dStkinf$. Again, for $\cU$ a $\infty$-topos, $\Fun(\Fun(\cX, \cX), \cU)$ is an $\infty$-topos as well.\\

Coming back to $\cX = \dStkinf$, which has a subobject classifier $\Omega$, one may ask whether this object itself is universal among all such subobject classifiers. Since we are dealing with $\infty$-topoi, we consider geometric morphisms between them. Let $\cE$ be an $\infty$-topos, and $\pi: \cE \rarr \dStkinf$ a geometric morphism, i.e. a left exact left adjoint. $\cE$ itself has a subobject classifier $G$ of its own, with $\pi G \in \cX$, such that:
\beq
\xymatrix{
S \ar@{>->}[d] \ar[r] &1 \ar[d] \\
\pi G \ar@{.>}[r] &\Omega
} \nonumber
\eeq
Now $\pi: \cE \rightleftarrows \dStkinf: i$ induces:
\beq
\pi_*: \Fun(\Delta^1, \cE) \rightleftarrows \Fun(\Delta^1, \dStkinf):i_* \nonumber
\eeq
From \cite{Lu1}, for $\cX$ an $\infty$-topos, $\cC$ any $\infty$-category, $\Fun(\cC, \cX)$ is an $\infty$-topos, so $\cO_{\cX}$ is an $\infty$-topos. $\pi_*$ is left exact, so we have again a geometric morphism:
\beq
\pi_*: \cO_{\cE} \rightleftarrows \cO_{\dStkinf} :i_* \nonumber
\eeq
It follows that $\pi_*$ carries monomorphisms in $\cE$ into monomorphisms in $\dStkinf$ (\cite{Lu1}). Further $\pi$ preserves pullbacks as a left adjoint, so $\pi_*$ preserves morphisms between monomorphisms, hence we have a well-defined morphism:
\beq
\pi_*: \cO_{\cE}^{(\mono)} \rarr \cO_{\dStkinf}^{(\mono)} \nonumber
\eeq
which is left exact, hence preserves final objects, in particular the subobject classifiers, so $\pi G \xrightarrow{\simeq} \Omega$. Now recall from \cite{Lu1} that for $K$ a simplicial set, one defines the homotopy category $hK$ of $K$ to be the homotopy category $h\gC[K]$ of the simplicial category $\gC[K]$. Finally for $K \in \SetD$, an object of $K$ is final if it is a final object of $hK$. In this sense $\pi G$ and $\Omega$ are isomorphic in $h \gC[\cO^{(\mono)}_{\dStkinf}] = h \cO^{(\mono)}_{\dStkinf}$, hence equivalent in $\dStkinf$.\\

Now the question is, what would be a universal such $\infty$-topos $\cE$. Recall that if one denotes by $\cR \cT op$ the category of $\infty$-topoi and right adjoints of geometric morphisms as morphisms between $\infty$-topoi, then $\Fun^*(\cS, \cX)$ is contractible for any $\infty$-topos $\cX$, which means exactly that $\cS$ is a final object of $\cR \cT op$. Consequently there is a geometric morphism $\dStkinf \rarr \cS$ with a left exact left adjoint, and by the same argument as above, if $\pi: \cS \rarr \dStkinf$ is a left exact left adjoint, it induces a left exact morphism $\cO^{(\mono)}_{\cS} \rarr \cO^{(\mono)}_{\dStkinf}$ mapping final object to final object, in a weak sense, which would make $\cS$ a universal $\infty$-topos in this regard, since $\cS$ is a final object of $\cR\cT op$, and consequently its subobject classifier the sought after universal subobject classifier.

\end{document}